\documentclass[a4paper, 11pt]{article}

\usepackage{algorithm}
\usepackage{algorithmic}
\usepackage{lineno,hyperref}
\usepackage{color}
\usepackage{amsfonts}
\usepackage{amsmath}
\usepackage{amssymb}
\usepackage{amsthm}
\usepackage{subeqnarray}
\usepackage{lscape,graphicx,url}
\usepackage{anysize}
\usepackage{proof}
\usepackage{fancyhdr}
\usepackage{latexsym}
\usepackage{array}
\usepackage{multirow}
\usepackage{textcomp}
\usepackage{mathtools}
\usepackage{proof}
\usepackage[margin=2cm]{geometry}
\usepackage{setspace}

\usepackage{tikz}
\usetikzlibrary{shapes,snakes}

\usepackage{environ}
\makeatletter
\newsavebox{\measure@tikzpicture}
\NewEnviron{scaletikzpicturetowidth}[1]{%
\def\tikz@width{#1}%
\def\tikzscale{1}\begin{lrbox}{\measure@tikzpicture}%
\BODY
\end{lrbox}%
\pgfmathparse{#1/\wd\measure@tikzpicture}%
\edef\tikzscale{\pgfmathresult}%
\BODY
}

\usepackage[none]{hyphenat}
\usepackage{paralist}
\usepackage{proof}
 \usepackage{placeins}
\usepackage{enumitem}
\usepackage{csquotes}
\usepackage{comment}
\usepackage{cite}
\usepackage{bm}
\usepackage{array}
\usepackage{rotating}

\usepackage{nicematrix}
\usepackage[misc]{ifsym}

\usepackage{etoolbox}
\patchcmd{\pprintMaketitle}
  {\fi\hrule}
  {\fi\ifvoid\extrainfobox\else\unvbox\extrainfobox\par\vskip10pt\fi\hrule}
  {}{}

\newsavebox\extrainfobox

\makeatletter
\def\thickhline{%
  \noalign{\ifnum0=`}\fi\hrule \@height \thickarrayrulewidth \futurelet
   \reserved@a\@xthickhline}
\def\@xthickhline{\ifx\reserved@a\thickhline
               \vskip\doublerulesep
               \vskip-\thickarrayrulewidth
             \fi
      \ifnum0=`{\fi}}
\makeatother

\newlength{\thickarrayrulewidth}
\newcolumntype{H}{>{\setbox0=\hbox\bgroup}c<{\egroup}@{}}

\newtheorem*{remark}{Theorem}

\newcommand{\GG}[1]{}

\begin{document}

\begin{center} 
{\Large An optimal algorithm for variable knockout problems} ~\\
~\\
{\Large J.E. Beasley} ~\\
~\\
Mathematics, Brunel University, Uxbridge UB8 3PH, UK ~\\
\Letter ~~john.beasley@brunel.ac.uk ~\\
~\\
February 2023,~Revised September 2023, October 2023\\
\end{center}

\begin{abstract}

We consider a class of problems related to variable knockout, where knockout means set a variable to zero. Given an
optimisation problem formulated as a zero-one integer program the question we consider in this paper
is what might be an appropriate set of variables to
knockout of the problem, in order that the optimal solution to the problem that remains
after variable knockout has a desired property. This property might be related to the optimal solution value
after knockout, 
or require the problem after knockout to be infeasible.
We present an algorithm for the optimal solution of this knockout problem. Computational results
are given for an illustrative example based upon shortest path interdiction
using publicly available shortest path test problems.

\end{abstract}


\section{Introduction} 

In this paper we consider problems related to variable knockout within the context of an optimisation problem formulated as a zero-one integer program. 

To illustrate the problem suppose that we have an optimisation problem involving $n$ zero-one variables $[x_i,~i=1,\ldots,n]$ and $m$ constraints where the optimisation problem is
\begin{equation}
\min  ~~\sum_{i=1}^n  c_i x_i
\label{eq1}
\end{equation}

\noindent subject to

\begin{equation}
 \sum_{j=1}^n  a_{ij} x_j \geq b_i \;\;\; i = 1, \ldots, m
\label{eq2}
\end{equation}
\begin{equation}
x_i \in \{ 0,1\} \;\;\; \;\;i=1,\ldots,n. 
\label{eq3}
\end{equation}

Then the class of problems considered in this paper relate to
what might be an appropriate set of variables to knockout of the problem, in order that the optimal solution to the problem that remains after variable knockout has a desired property. Here by variable knockout we mean explicitly set any variable knocked out to zero. So if $D$ is the set of variables that are knocked out we have the constraint:
\begin{equation}
  x_i = 0 \;\;\; \forall i \in D.
\label{eq4}
\end{equation}

So the optimisation problem with the knocked out variables is now optimise Equation~(\ref{eq1}) subject to Equations~(\ref{eq2})-(\ref{eq4}). We wish to 
choose $D$ such that the optimal solution to this problem has a desired property. This property might be:
\begin{enumerate}
\item \emph{\textbf{P1}}: the optimised (minimal) objective function value after knockout is at least (so greater than or equal to) a particular value; or

\item\emph{\textbf{P2}}: the optimisation problem that remains after knockout is infeasible (i.e.~one or more constraints cannot be satisfied).
\end{enumerate}

\noindent To formulate the variable knockout problem 
introduce variables $[\alpha_i,~i=1,\ldots,n]$ where $\alpha_i=1$ if variable $i$ is knocked out, zero otherwise. Then we have
\begin{equation}
x_i + \alpha_i \leq 1 \;\;\; i=1,\ldots,n 
\label{eq11}
\end{equation}
\begin{equation}
 \alpha_i  \in \{ 0,1\} \;\;\; i=1,\ldots,n. 
\label{eq11a}
\end{equation}

Equation~(\ref{eq11}) ensures that if $\alpha_i=1$ then $x_i=0$, so the variable is knocked out. If $\alpha_i=0$ then 
Equation~(\ref{eq11}) has no effect on the value of $x_i$ adopted, since it can be either zero or one.

\sloppy Clearly for either of the two properties P1 and P2 given above
there might be many sets of knocked out variables $D$ such that the resulting optimisation problem, optimise Equation~(\ref{eq1}) subject to Equations~(\ref{eq2})-(\ref{eq4}), has the required property. However, since we are in an  optimisation context, it would be natural to assign a value $d_i$ to each variable $i$ that is knocked out and  consider the problem of choosing the set $D$ such that it  minimises $\sum_{i \in D } d_i \alpha_i$. For example setting $d_i=1~i=1,\ldots,n$ would correspond to choosing the minimal set of variables $D$ to knockout such that the resulting optimisation problem has the desired property.

Considering the first property P1 above suppose that we wish to knockout variables such that the optimal solution to the problem that remains after knockout has value of at least $C^*$.
The variable knockout problem considered in this paper  associated with P1 can then be stated as
\begin{equation}
\min  ~~\sum_{i=1}^n  d_i \alpha_i
\label{jebex1}
\end{equation}
\noindent subject to 
\begin{equation}
\bigg\{\mbox{optimise~Equation~(\ref{eq1})~subject~to~Equations~(\ref{eq2}),(\ref{eq3}),(\ref{eq11}),(\ref{eq11a})} \bigg\}~\mbox{has~value~$\geq~C^*$}. 
\label{eqjebadd}
\end{equation}

Considering the second property P2 above  suppose that we wish to knockout variables such that  the problem that remains after knockout is infeasible.
The variable knockout problem considered in this paper  associated with P2 can then be stated as
\begin{equation}
\min  ~~\sum_{i=1}^n  d_i \alpha_i
\label{jebex1p2}
\end{equation}
\noindent subject to
\begin{equation}
\bigg\{\mbox{optimise~Equation~(\ref{eq1})~subject~to~Equations~(\ref{eq2}),(\ref{eq3}),(\ref{eq11}),(\ref{eq11a})} \bigg\}~\mbox{is infeasible.} 
\label{eqjebaddp2}
\end{equation}

The motivation for the author to develop the algorithms given in this paper grew out of his work in systems biology 
with regard to metabolic pathways~\cite{bf15,bf09,bf09a,bf07}. These are networks of biochemical reactions in organisms 
and in that context knockout refers to neutralising reactions to affect the overall system (e.g.~to enhance/reduce the production of certain chemicals).

The structure of this paper is as follows. 
In Section~\ref{lit} we discuss work in the literature that is relevant to the variable knockout problem considered in this paper. This works deals with two topics: bilevel optimisation and interdiction problems. We also discuss what we believe is the contribution of this paper to the literature. 
In Section~\ref{algorithm} we present our algorithm for the optimal solution of the variable knockout problem and prove its correctness.
We also 
 indicate
how our algorithm can be adapted when the number of variables knocked out is specified (i.e.~when we have a cardinality constraint). 
In Section~\ref{Results} computational results are given for an illustrative shortest path example involving arc knockout. We also present results for shortest path cardinality constrained knockout. 
Finally in Section~\ref{Sec:Conclusions} we present our conclusions.

\section{Literature: bilevel optimisation and interdiction problems}
\label{lit}

\subsection{Bilevel optimisation}
In bilevel optimisation~\cite{schmidt23} we have an upper level  optimisation problem, often called the leader's problem, and a lower level
optimisation problem, often called the follower's problem, where the lower-level problem involves the decisions made at the upper level.

\sloppy The zero-one integer program, Equations~(\ref{jebex1}),(\ref{eqjebadd}), given above is an integer 
(zero-one) bilevel optimisation problem, 
with the $[\alpha_i]$ variables involved in both the upper  level optimisation (Equation~(\ref{jebex1})) as well as in the 
lower level optimisation (Equation~(\ref{eqjebadd}), since that involves Equation~(\ref{eq11})). 
Bilevel optimisation has been a subject of interest for many years, e.g.~\cite{denegre09, Colson2007, Dempe2003, Kalashnikov2015, sinha18, Vicente94} and for recent papers that survey bilevel optimisation see~\cite{beck23, dios23, kleinert21}.

Integer bilevel optimisation
 problems, so bilevel problems in which all variables are integer, are generally challenging to solve.
However in this instance our bilevel program has a special structure, namely a 
constraint on the value of the objective function, Equation~(\ref{eqjebadd}), and this enables us to present an algorithm for its optimal solution.
As best as we are aware knockout problems of the type considered here have not been significantly explored in the literature.

\subsection{Interdiction problems}
Interdiction problems are a special class of bilevel problems~\cite{kleinert21, smith20}. In problems of this type the typical interpretation is that the leader is 
attempting, by 
their actions, to block or inhibit the activities of the follower.  The classical example here is that the leader is taking action to 
remove arcs/nodes from a network so as to affect the shortest path that
 the follower can use between a specified origin and destination~\cite{Israeli2002}.
This network interdiction problem has been much studied, e.g.~see~\cite{Israeli2002, Rocco2010, Bayrak2008, smith20, kleinert21}. However, interdiction problems in which the follower's decision problem is not a shortest path related problem also exist, for example the follower's problem can be a knapsack problem, a maximum clique problem, a maximum matching problem;
 see~\cite{smith20}.

The optimisation problems associated with P1 and P2 above are interdiction problems since setting to one
a knockout variable associated with the leader (upper level optimisation)  forces (via Equation~(\ref{eq11})) the 
associated $x$ variable (in the follower lower level optimisation) to 
assume the value zero. 

\subsection{Contribution}
We believe that the primary contribution of this paper to the literature is:
\begin{itemize}
\item to introduce the knockout problem, where (in bilevel optimisation terms) we have a constraint upon the optimal solution of the follower's
decision problem
\item to present an algorithm for the optimal solution of the knockout problem
\item to demonstrate the computational performance of this algorithm using publicly available shortest path test problems.
\end{itemize}
As best as we are aware knockout problems of the type considered here have not been significantly explored in the literature.
A secondary contribution is
\begin{itemize}
\item  to demonstrate that the algorithm given for the optimal solution of the knockout problem can be easily extended to cardinality constrained knockout
\item to demonstrate the computational performance of this algorithm for cardinality constrained knockout using publicly available shortest path test problems.
\end{itemize}

\section{Solution algorithms} \label{algorithm}

In this section we first present our algorithm for solving the variable knockout problem optimally and prove its correctness.
We then indicate
how our algorithm can be adapted when the number of variables knocked out is specified (i.e.~when we have a cardinality constraint).

\subsection{Optimal algorithm}
The pseudocode for our algorithm for optimally solving the knockout problem is shown in Algorithm~\ref{alg1}. This pseudocode deals with both desired properties (P1 and P2, as discussed above).


Let $F(s)$ be the set of non-zero variables $[i~|~x_i=1~i=1\ldots,n]$ in feasible solution $s$ ($s=1,\ldots,S$).
In Algorithm~\ref{alg1} we start by first solving the problem being considered and initialising the number of feasible solutions ($S$) to one, with the corresponding solution being $F(1)$. 

Given a solution $F(1)$ then it is clear that in order to avoid this solution being active we need to knockout at least one of the non-zero variables in $F(1)$. So we need to add an appropriate knockout constraint $\sum_{i \in F(1)} \alpha_i \geq 1$. 

As long as we do not have a solution with the required property 
we continue solving, 
but with knockout constraints added for all previously identified solutions. Each new solution found is included as $F(S)$. Note here that as 
we are continually adding a new knockout constraint the value of the objective function value, Equation~(\ref{eq1}), can never decrease at 
each iteration, rather it must increase (or remain the same).

Once we have a solution with the required property then we find the optimal knockout set. The problem to be solved here involves the
 original problem constraints if the desired property is P1, since in that case we need a feasible solution. If the desired property is P2, so the 
original problem is infeasible after knockout, then we do not need to include the original problem constraints in the optimisation.

Note here that the optimisation problem associated with P2 to determine the optimal knockout set, namely 
optimise $\sum_{i=1}^n d_i\alpha_i$ subject to Equation~(\ref{eq11a}) and $\sum_{i \in F(s)} \alpha_i \geq 1~s=1,\ldots,S$
is a set covering problem involving $n$ variables and $S$ constraints. The set covering problem
has been extensively studied in the literature and can be solved, either  optimally or heuristically, for very large
instances~\cite{beasley92, lan07, caprara00, caprara99, beasley96, beasley87,reyes21,naji10}.

\begin{algorithm}[!htb]
\phantom{a}
\caption{Pseudocode for the knockout algorithm}
\begin{algorithmic}[1]
\STATE Optimise Equation~(\ref{eq1}) subject to Equations~(\ref{eq2}),(\ref{eq3})   \hfill  \COMMENT{Solve the problem under consideration}\\
\STATE $S \leftarrow 1$ \\
\STATE $F(S) \leftarrow$ non-zero variables in the current solution \hfill  \COMMENT{Record the solution}\\
\STATE $t \leftarrow 1$ \hfill  \COMMENT{Initialise the iteration counter}\\
\WHILE{solution does not have the desired property} 
\STATE $t \leftarrow t +1$ \hfill  \COMMENT{Update the iteration counter}\\
\STATE Optimise Equation~(\ref{eq1}) subject to Equations~(\ref{eq2}),(\ref{eq3}),(\ref{eq11}),(\ref{eq11a}) and  $\sum_{i \in F(s)} \alpha_i \geq 1~s=1,\ldots,S$   \\
   \hfill \COMMENT{Solve the problem with the knockout constraints added} \\
\STATE If the desired property is P1 and there is no feasible solution in line 7 then it is not possible to achieve the desired property so terminate \\
\hfill  \COMMENT{Check for infeasibility when considering P1}\\
\STATE  $S \leftarrow S+1 $  \\
\STATE $F(S) \leftarrow$ non-zero variables in the current solution    \hfill 
\COMMENT{Record the solution} \\
\ENDWHILE
\STATE If the desired property is P1 then: Optimise Equation~(\ref{jebex1}) subject to Equations~(\ref{eq2}),(\ref{eq3}),(\ref{eq11}),(\ref{eq11a}) and $\sum_{i \in F(s)} \alpha_i \geq 1~s=1,\ldots,S$   \\ \hfill \COMMENT{Find the optimal knockout set}\\
\STATE If the desired property is P2 then: Optimise Equation~(\ref{jebex1p2}) subject to Equation~(\ref{eq11a}) and $\sum_{i \in F(s)} \alpha_i \geq 1~s=1,\ldots,S$   \\ \hfill \COMMENT{Find the optimal knockout set}\\
\end{algorithmic}
  \label{alg1}
\medskip
\end{algorithm}

\begin{remark}
Algorithm~\ref{alg1} finds the optimal solution to the variable knockout problem
associated with P1 
 (Equations~(\ref{jebex1}),(\ref{eqjebadd})) 
in a finite number of iterations.
\end{remark}

\noindent A detailed and comprehensive proof of this theorem is given below.

\begin{proof} 
$ $\newline
\begin{compactitem}
\vspace{-5mm} 

\item \emph{\textbf{Line 1}} We solve the original problem. Assume, without significant loss of generality, that this solution $F(1)$ has value $< C^*$, so does not have the desired property, Equation (\ref{eqjebadd}).
Clearly this solution must be eliminated to satisfy Equation (\ref{eqjebadd}) and this is accomplished via use of the constraint $\sum_{i \in F(1)} \alpha_i \geq 1$ in line 7.

\item \emph{\textbf{Lines 5-11}} In the $\mbox{while} \ldots \mbox{end~while}$ loop we repetitively solve, each time adding to $F(S)$ in line 10 a solution that must be eliminated in order to satisfy Equation (\ref{eqjebadd}). Hence the use of the constraint $\sum_{i \in F(s)} \alpha_i \geq 1~s=1,\ldots,S $ in line 7. Here  (assuming we have a feasible solution) we exit the loop as soon as we detect at line 7 a solution that does not have the desired property (so we do not go to lines 8-10 and then exit).

\item \emph{\textbf{Line 8}} If the optimisation problem in line 7 is infeasible, and we are dealing with P1, then we terminate since it is impossible to find a feasible solution with the desired P1 property (so with value $\geq C^*$).
\item \emph{\textbf{Line 12}} When we exit the $\mbox{while} \ldots \mbox{end~while}$ loop it must be the case that the current set $F(S)$ of added solutions are such that it is 
 \emph{\textbf{impossible}}
 to find a solution which will violate Equation (\ref{eqjebadd}), so with a solution value $<C^*$, since if such a solution exists it would have been
detected at line 7 in the previous iteration of the $\mbox{while} \ldots \mbox{end~while}$ loop. \emph{\textbf{In other words since we progressively find feasible solutions with
non-decreasing solution values it must be the case for P1 that when we exit the $\mbox{while} \ldots \mbox{end~while}$ loop the current solution 
 \boldmath
$F(S)$ 
 \unboldmath 
has solution value 
 \boldmath
$\geq C^*$.}}
 \unboldmath 
\item \emph{\textbf{Line 12}} The final optimisation takes the solution set $F(s)~s=1,\ldots,S$ and finds an optimal knockout set. This knockout set will still have the desired property (so with solution value $\geq C^*$) because as we have exited the 
$\mbox{while} \ldots \mbox{end~while}$ loop
in the final iteration we know that are no possible knockout sets that violate Equation (\ref{eqjebadd}), so with solution value $< C^*$. 
\item Hence we have proved that Algorithm~\ref{alg1} finds the optimal solution to the variable knockout problem
associated with P1 
 (Equations~(\ref{jebex1}),(\ref{eqjebadd})).

\end{compactitem}
 \noindent The number of iterations ($t$) required must be finite since at each iteration in the  $\mbox{while} \ldots \mbox{end~while}$ loop we find a new feasible solution and the entire 
set of feasible solutions are a subset of the  $2^n$ possible combinations of values for the zero-one variables $[x_1,x_2,\ldots,x_n]$.
\end{proof}
\noindent Note here that, as far as we are aware, the  algorithm  given above for the optimal solution of the variable knockout problem  has not been presented previously in the literature. Given the proof above for P1 it is trivial to see that the algorithm presented  is also optimal for solving with regard to P2.
Note also here that:
\begin{compactitem}
\item We do not need to assume that the solution $F(S)$ in Algorithm~\ref{alg1} is uniquely defined, if multiple optimal solutions exist choice of any one of them will suffice.

\item
The knockout constraints added 
in Algorithm~\ref{alg1} are of  the form $\sum_{i \in F(s)} \alpha_i \geq 1~s=1,\ldots,S$ and can be regarded as \emph{\textbf{covering constraints}} (since constraints of this type appear in set covering problems). Such constraints 
 have also appeared in other papers in the literature. For example in~\cite{Israeli2002}  in a Benders decomposition
 approach to shortest path interdiction 
(where they are referred  to  as \enquote{supervalid inequalities}), 
and where they are lifted (strengthened) by making use of lower bound information. In \cite{lozano17} covering constraints are added associated with fortification of assets before interdiction,
see also\cite{leitner23}.  Wei and Walteros~\cite{wei22} give a number of examples of covering constraints for different interdiction problems.

\item Our knockout constraints are related to \enquote{no-good cuts}. For example consider $F(1)$, the set of 
non-zero variables in the first solution. 
The standard no-good cut to eliminate this solution is $\sum_{i \notin F(1)} x_i + \sum _{i \in F(1)}(1- x_i) \geq 1$.
The knockout constraint is $\sum_{i \in F(1)} \alpha_i \geq 1$, which from Equation~(\ref{eq11}) is equivalent to forcing at least one of the $[x_i~|~i \in F(1)]$ to zero, 
i.e.~equivalent to $ \sum _{i \in F(1)}(1- x_i) \geq 1$.
Hence the knockout constraint is a no-good cut tightened by removing the $\sum_{i \notin F(1)} x_i$ term.

\item The optimisation problem (involving the original objective, Equation~(\ref{eq1})) solved in Algorithm~\ref{alg1} involves 
\emph{\textbf{simultaneously}} deciding values for the original problem variables $[x_i]$ as well as the knockout variables $[\alpha_i]$. 

\item It is simple to show by means of an example that our knockout algorithm does not involve complete enumeration of all possible feasible solutions lacking the desired property. This is because the addition of a single knockout constraint can remove from consideration more than one feasible solution.
\end{compactitem}

\noindent Clearly the computational effectiveness of our knockout algorithm in identifying the optimal set of knockout variables will depend upon both  the underlying problem under consideration and the required property. For example it is clear that
identifying knockout variables to render the problem infeasible (so a problem of type P2) would be more challenging than identifying knockout variables that simply raise the objective function value slightly from the minimal value achieved with no knockout (so a problem of type P1).


\subsection{Cardinality constrained knockout}
Suppose that we wish to constrain the number of variables knocked out, i.e.~impose the cardinality constraint
\begin{equation}
\sum_{i=1}^n \alpha_i = K. 
\label{eqcard}
\end{equation}
\noindent Here we need $K \leq n$ since clearly we cannot knockout more variables than we have in the problem. 

In this case \textbf{\emph{an obvious objective is 
 to choose the $K$ best variables to knockout so as 
to maximise the (minimal) value of original problem after elimination of the knocked out variables}}. So this problem is
\begin{equation}
\small
\mbox{maximise the solution value of}~\bigg\{\mbox{optimise~Equation~(\ref{eq1})~subject~to~Equations~(\ref{eq2}),(\ref{eq3}),(\ref{eq11}),(\ref{eq11a}),(\ref{eqcard})} \bigg\}
\label{eqcardopt}
\end{equation}
\normalsize

In order to solve this problem, Equation~(\ref{eqcardopt}), we simply modify Algorithm~\ref{alg1}, as shown in Algorithm~\ref{alg2}. 
The logic underlying Algorithm~\ref{alg2} is similar to that underlying 
Algorithm~\ref{alg1}. In order to maximise the (minimal) value of the original problem after variable knockout we need to successively eliminate (knockout) solutions $F(1)$, $F(2)$, etc; where we can only knockout $K$ variables in total (Equation~(\ref{eqcard})).

Algorithm~\ref{alg2} is predominantly the same as Algorithm~\ref{alg1} except that we add Equation~(\ref{eqcard}) to the optimisation
and repeat the solution process 
until the problem is infeasible. 
We then know that adding knockout constraints for the  $S$ solutions found, in conjunction with the cardinality constraint
(Equation~(\ref{eqcard})), renders the problem infeasible. Hence the maximal value that can be achieved for the (minimised) objective can be found by solving the problem with 
the first   $(S-1)$ solutions found knocked out. 

The proof that Algorithm~\ref{alg2} finds the optimal solution  to the cardinality constrained variable knockout problem
is similar to the proof for Algorithm~\ref{alg1} given above and so, for space reasons, is not given here.

On a technical note here Algorithm~\ref{alg2} maximises the (minimal) value of a \textbf{\emph{feasible solution}} to the original problem when $K$ variables are knocked out. If a feasible solution exists with $K$ variables knocked out then it is possible that the original problem could also be rendered infeasible by judicious choice of $K$ (or fewer) variables to knockout, but this would not 
be detected by Algorithm~\ref{alg2}. 

\begin{algorithm}[!htb]
\phantom{a}
\caption{Pseudocode for cardinality constrained knockout}
\begin{algorithmic}[1]
\STATE Optimise Equation~(\ref{eq1}) subject to Equations~(\ref{eq2}),(\ref{eq3})   \hfill  \COMMENT{Solve the problem under consideration}\\
\STATE $S \leftarrow 1$ \\
\STATE $F(S) \leftarrow$ non-zero variables in the current solution \hfill  \COMMENT{Record the solution}\\
\STATE $t \leftarrow 1$ \hfill  \COMMENT{Initialise the iteration counter}\\
\WHILE{problem is not infeasible} 
\STATE $t \leftarrow t +1$ \hfill  \COMMENT{Update the iteration counter}\\
\STATE Optimise Equation~(\ref{eq1}) subject to Equations~(\ref{eq2}),(\ref{eq3}),(\ref{eq11}),(\ref{eq11a}),(\ref{eqcard})
 and  $\sum_{i \in F(s)} \alpha_i \geq 1~s=1,\ldots,S$   \\
   \hfill \COMMENT{Solve the problem with the knockout constraints added} \\
\STATE  $S \leftarrow S+1 $  \\
\STATE $F(S) \leftarrow$ non-zero variables in the current solution    \hfill 
\COMMENT{Record the solution} \\
\ENDWHILE
\STATE  $S \leftarrow S-1 $  
\hfill  \COMMENT{Reduce $S$ by one} \\
\STATE Optimise Equation~(\ref{eq1}) subject to Equations~(\ref{eq2}),(\ref{eq3}),(\ref{eq11}),(\ref{eq11a}),(\ref{eqcard})
 and  $\sum_{i \in F(s)} \alpha_i \geq 1~s=1,\ldots,S$  \\ \hfill \COMMENT{Find the optimal knockout set and objective function value}\\
\end{algorithmic}
  \label{alg2}
\medskip
\end{algorithm}

\section{Illustrative example} \label{Results}

As mentioned previously above the motivation for the author to develop the algorithms given above grew out of his work in systems biology 
with regard to metabolic pathways~\cite{bf15,bf09,bf09a,bf07}. These are networks of biochemical reactions in organisms 
and in that context knockout refers to neutralising reactions to affect the overall system (e.g.~to enhance/reduce the production of certain chemicals).

Clearly few Operational Researchers will have any familiarity with systems biology so
to illustrate computationally the optimal knockout algorithm presented in this paper we consider the problem of finding
 the minimal number of arcs to knockout from a directed network such that, after knockout,  the shortest path from an origin 
node to a destination node is of length at least a specified value. Problems of this type, so shortest path network interdiction, have been considered 
previously in the literature,
e.g.~see~\cite{Israeli2002, Rocco2010, Bayrak2008, smith20, kleinert21}.

\emph{\textbf{It is very important to emphasise here that we are using this shortest path network interdiction problem purely to 
illustrate computationally our general optimal knockout algorithm.
We are not suggesting that the general knockout algorithm presented in this paper is a computationally competitive approach
for shortest path network interdiction. Other approaches for shortest path network interdiction typically make 
use  of problem-specific knowledge.
By contrast our general knockout algorithm makes no  use of problem-specific knowledge.
Here we are simply illustrating the application of our general knockout algorithm using shortest path network 
interdiction as an example problem.}}

\sloppy For our illustrative purposes we made use of 
a small number of directed network instances given (albeit in a different context) in~\cite{beasley89}, which have the advantage for future workers studying knockout algorithms that these instances are publicly available from OR-Library~\cite{beasley90}, see
\href{http://people.brunel.ac.uk/~mastjjb/jeb/info.html}{http://people.brunel.ac.uk/$\sim$mastjjb/jeb/info.html}.
The computational results presented below 
(Windows  2.50GHz pc, Intel i5-2400S processor, 6GB memory)
are for our optimal knockout algorithms (Algorithm~\ref{alg1}
and Algorithm~\ref{alg2})
 as coded in FORTRAN using SCIP (Solving Constraint Integer Programs)
\cite{Achterberg2009,scip} as the 
integer solver. 

Table~\ref{table2} shows the results obtained for Algorithm~\ref{alg1}. In that table we show the number of nodes and arcs in the instances considered, as well as the length of the shortest path from the origin node to the destination node. We considered the problem of  knocking out the minimal  number of arcs such that, after knockout,  the length of the remaining shortest path from the origin node to the destination node was at least a multiplier $\gamma$ of the length of the original (unknocked out) shortest path. Table~\ref{table2} shows results for $\gamma{=}1.5$ and $\gamma{=}2$.

In that table we can see that for the first problem with 100 nodes and 955 arcs the shortest path (from node 1 to node 100) is of length 80. The minimal number of arcs to knockout in order that the shortest path is of length at least $1.5(80){=}120$ is 2 and our algorithm required (in total) 17.4 seconds to prove this, with $S{=}21$ solutions being found. The minimal number of arcs to knockout in order that the shortest path is of length at least $2(80){=}160$ is 4, requiring (in total)  783.3  seconds, with $S{=}281$ solutions being found.

\begin{table}[!htb]
{\scriptsize
\begin{tabular}{cccccccccc}
\hline
Number  & Number & Shortest &  \multicolumn{3}{c}{Multiplier $\gamma{=}1.5$}
&  & \multicolumn{3}{c}{Multiplier $\gamma{=}2$} \\
\cline{4-6} \cline{8-10}
of nodes & of arcs & path  & Number of & Minimal number  & Total time & &  Number of & Minimal number  & Total time \\
 & & length & solutions  & of arcs  & (secs) &  & solutions  & of arcs  & (secs) \\
& & & ($S$) & knocked out & & & ($S$) & knocked out  \\
\hline

100	&	955	&	80	&	21	&	2	&	17.4	&	&	281	&	4	&	783.3	  \\	   
	&	990	&	79	&	10	&	1	&	5.8	&	&	135	&	4	&	310.9	  \\	   
200	&	2040	&	230	&	30	&	4	&	113.9	&	&	852	&	7	&	10532.4	  \\	   
	&	2080	&	200	&	16	&	3	&	42.5	&	&	322	&	6	&	2305.8	  \\	   
500	&	4858	&	455	&	2	&	1	&	44.8	&	&	28	&	3	&	445.3	  \\	   
	&	4847	&	611	&	18	&	5	&	284.9	&	&	1325	&	8	&	53061.9	  \\	 

\hline 	   
								 
\end{tabular}
\caption{Shortest path results}\label{table2}
}
\end{table}

Table~\ref{table3} shows, for the same problems as considered in Table~\ref{table2}, the results for cardinality constrained knockout (Algorithm~\ref{alg2}) for $K=1,2,3,4,5,10$. To illustrate this table consider the first problem  with 100 nodes, 955 arcs and shortest path length 80. For $K=1$, i.e.~a single arc knocked out so as to maximise the length of the shortest path in the network that remains after knockout, the shortest path length increases to 110, requiring $S=3$ solutions and 3.2 seconds in total.  For $K=2$, i.e.~two arcs knocked out so as to maximise the length of the shortest path in the network that remains after knockout, the shortest path length increases to 139, requiring $S=6$ solutions and 5.2 seconds in total. 
Examining Table~\ref{table3} we can see that, as we would expect, for any particular problem as $K$ increases the number of solutions $S$ also increases.

\begin{table}[!htb]

\begin{tabular}{llcccccc}
\hline

		\multicolumn{2}{l}{Number of nodes}	&	100	&	100	&	200	&	200	&	500	&	500			  \\		   
		\multicolumn{2}{l}{Number of arcs}	&	955	&	990	&	2040	&	2080	&	4858	&	4847			  \\		   
		\multicolumn{2}{l}{Shortest path length}	&	80	&	79	&	230	&	200	&	455	&	611			  \\		   
 
\hline																		   
$K=1$	&	Number of solutions ($S$)	&	3	&	2	&	2	&	3	&	3	&	1			  \\	   
	&	Shortest path length	&	110	&	119	&	260	&	258	&	779	&	689			  \\	   
	&	Total time (secs)	&	3.2	&	2.6	&	10.4	&	14.4	&	76.7	&	51.7			  \\	   
$K=2$	&	Number of solutions ($S$)	&	6	&	5	&	6	&	4	&	9	&	2			  \\	   
	&	Shortest path length	&	139	&	122	&	308	&	266	&	906	&	715			  \\	   
	&	Total time (secs)	&	5.2	&	4.7	&	20.7	&	15.5	&	169.8	&	69.3			  \\	   
$K=3$	&	Number of solutions ($S$)	&	12	&	10	&	11	&	7	&	15	&	6			  \\	   
	&	Shortest path length	&	142	&	154	&	321	&	317	&	913	&	838			  \\	   
	&	Total time (secs)	&	8.8	&	7.6	&	34.0	&	23.4	&	275.4	&	118.7			  \\	   
$K=4$	&	Number of solutions ($S$)	&	22	&	18	&	20	&	11	&	24	&	9			  \\	   
	&	Shortest path length	&	185	&	212	&	360	&	334	&	986	&	866			  \\	   
	&	Total time (secs)	&	19.9	&	13.8	&	61.6	&	31.4	&	473.4	&	168.5			  \\	   
$K=5$	&	Number of solutions ($S$)	&	40	&	30	&	34	&	14	&	31	&	16			  \\	   
	&	Shortest path length	&	209	&	232	&	418	&	339	&	1070	&	979			  \\	   
	&	Total time (secs)	&	60.6	&	34.7	&	164.9	&	41.9	&	671.6	&	294.0			  \\	   
$K=10$	&	Number of solutions ($S$)	&	255	&	160	&	267	&	136	&	260	&	144			  \\	   
	&	Shortest path length	&	263	&	275	&	619	&	498	&	1334	&	1389			  \\	   
	&	Total time (secs)	&	833.7	&	531.5	&	2908.8	&	913.8	&	11130.7	&	3651.3			  \\	   
\hline

\end{tabular}
\caption{Cardinality constrained knockout results}\label{table3}

\end{table}

\section{Conclusions} \label{Sec:Conclusions}

In this paper we have considered a class of problems related to variable knockout. In problems of this type we need to 
decide an appropriate set of variables to knockout of the problem, in order that the optimal solution to the problem that remains after variable knockout has a desired property.

We presented an algorithm for the optimal solution of the problem. We also  illustrated 
how our algorithm could be adapted when the number of variables knocked out is specified (i.e.~when we have a cardinality constraint).

To illustrate our knockout approach we used publicly available test problems.
Computational results were given for the problem of finding the minimal number of arcs to knockout from a directed network
such that, after knockout,  the shortest path from an origin node to a destination node was of length at least a specified value. 
We also presented results for shortest path cardinality constrained knockout.

\appendix
\section*{Appendix: Example}
\label{jebappa}
In this Appendix we illustrate Algorithm~\ref{alg1} and  Algorithm~\ref{alg2} using the example directed graph with five nodes shown in Figure~\ref{fig1}, where 
the underlying minimisation problem is to find the shortest path from node 1 to node 5. In this problem variable knockout corresponds to eliminating an arc from the problem.

For readers unfamiliar with the standard formulation of the shortest path problem it is as described here.
Let $A$ be the set of arcs in Figure~\ref{fig1}, so $A= [(1,2), (1,3), (2,5), (3,4), (3,5), (4,5)]$ with $c_{ij}$ being the cost of arc $(i,j)$ as shown by the number next to each arc 
in  Figure~\ref{fig1}.
The associated decision variable is $x_{ij}=1$ if arc $(i,j) \in A$ is used, 0 otherwise. Then Equation~(\ref{eq1}) is
\begin{equation*}
\min  ~~\sum_{(i,j) \in A} c_{ij}x_{ij}
\end{equation*}
\noindent and the constraints associated with Equation~(\ref{eq2}) are
\begin{equation*}
x_{12} + x_{13} = x_{25} + x_{35} + x_{45} = 1
\end{equation*}
\begin{equation*}
x_{12} = x_{25}
\end{equation*}
\begin{equation*}
x_{13} = x_{34} + x_{35}
\end{equation*}
\begin{equation*}
x_{34} = x_{45}.
\end{equation*}
\noindent The first constraint here ensures just one arc comes out of node 1 and one arc goes into node 5. The other constraints ensure that the number of arcs into nodes 2, 3 and 4 respectively are equal to the number of arcs out.

\begin{figure}[!htb]
\centering	  
\begin{scaletikzpicturetowidth}{\textwidth*0.70} 
\begin{tikzpicture}[scale=\tikzscale,auto=center,every node/.style={circle, draw, ultra thick}]

 \node (a1) at (0,1) {1};  
  \node (a2) at (1,2)  {2}; 
  \node (a3) at (1,0)  {3};  
  \node (a4) at (2,0) {4};  
  \node (a5) at (3,1)  {5};  
  

     \draw[->, ultra thick] (a1) -- (a2) node [draw=none, left, above, midway] {1};

    \draw[->, ultra thick] (a1) -- (a3) node [draw=none, left, below, midway] {2};  

    \draw[->, ultra thick] (a2) -- (a5) node [draw=none, left, above, midway] {1};  

     \draw[->, ultra thick](a3) -1 (a5) node [draw=none, left, above, midway] {1};


    \draw[->, ultra thick] (a3) -1 (a4) node [draw=none, left, below, midway] {1};
    \draw[->, ultra thick] (a4) -1 (a5) node [draw=none, left, below, midway] {1};
  
\end{tikzpicture}  
\end{scaletikzpicturetowidth}
\caption{Example problem}
\label{fig1}
\end{figure}

To illustrate Algorithm~\ref{alg1} suppose we consider $d_i=1~\forall i$ (so Equation~(\ref{jebex1}) becomes $\sum_{(i,j) \in A} \alpha_{ij}$) and $C^*=4$. 
So here we have a problem of type P1 and are seeking the minimal number of arcs to knockout 
such that the value of the shortest path after their elimination is at least 4.
\begin{compactitem}
\item In line 1 of Algorithm~\ref{alg1} we first solve the original shortest path problem with the path chosen being 1-2-5 of cost 2.
\item Then in  line 7 of
Algorithm~\ref{alg1} 
we again find the shortest path, but now where we have added the constraints given in Equation~(\ref{eq11}), which here are 
$x_{ij} + \alpha_{ij} \leq 1~\forall (i,j) \in A$,
together with
$\alpha_{12} + \alpha_{25} \geq 1$
and the integrality constraint Equation~(\ref{eq11a}).

\item  The path chosen is now 1-3-5 of cost 3, where there are many optimal solutions involving the $\alpha$ variables, one of which is $\alpha_{12}  =1$. The new constraint to be added is now $\alpha_{13} + \alpha_{35} \geq 1$.
\item Relooping to line 7 we again find the shortest path. The path chosen is 1-3-4-5 of cost 4 where $\alpha_{35} =1$ and either $\alpha_{12}$ or $\alpha_{25}$ is also one.
Since the solution is of value 4, so $\geq C^*$ we are done and go to line 12
with $F(1)=[x_{12},x_{25}], F(2)=[x_{13},x_{35}], S=2$.
\item In line 12 
of Algorithm~\ref{alg1} 
the objective is to minimise $\sum_{(i,j) \in A} \alpha_{ij}$ subject to the shortest path constraints (to ensure we have a feasible path) and
 $\alpha_{12} + \alpha_{25} \geq 1$ and $\alpha_{13} + \alpha_{35} \geq 1$. The solution is of value two.
 There are two optimal solutions here, with  $\alpha_{35}=1$ and  either $\alpha_{12}$ or $\alpha_{25}$  is also  one, with the associated path being  1-3-4-5 of cost 4. 
\end{compactitem}
\noindent So  the minimal number of arcs to knockout 
such that the value of the shortest path after their elimination is at least 4 is two, namely arc (3,5) and either arc (1,2) or arc (2,5). \\

 To illustrate Algorithm~\ref{alg2} suppose we consider $K=2$, so we are seeking the two best arcs to knockout 
so as to maximise the value of the shortest path after their elimination.

\begin{compactitem}
\item In line 1 of Algorithm~\ref{alg2} we first solve the original shortest path problem with the path chosen being 1-2-5 of cost 2.
\item Then in  line 7 of Algorithm~\ref{alg2} we again find the shortest path, but now where we have added the constraints given in 
Equation~(\ref{eq11}), which here are 
$x_{ij} + \alpha_{ij} \leq 1~\forall (i,j) \in A$,
and Equation~(\ref{eqcard}), which is $\sum_{(i,j) \in A} \alpha_{ij} =2$,
together with
$\alpha_{12} + \alpha_{25} \geq 1$ and the integrality constraint Equation~(\ref{eq11a}).
\item The path chosen is now 1-3-5 of cost 3, where there are many optimal solutions involving the $\alpha$ variables, one of which is $\alpha_{12} = \alpha_{25} =1$. The new constraint to be added is now $\alpha_{13} + \alpha_{35} \geq 1$.
\item Relooping to line 7 we again find the shortest path. The path chosen is 1-3-4-5 of cost 4 where $\alpha_{35} =1$ and either $\alpha_{12}$ or $\alpha_{25}$ is also one.
The new constraint to be added is now $\alpha_{13} + \alpha_{34}  + \alpha_{45} \geq 1$.
\item Relooping to line 7 the problem is infeasible so we go to line 11 with  $F(1)=[x_{12},x_{25}], F(2)=[x_{13},x_{35}], F(3)=[x_{13},x_{34},x_{45}],
S=3$.
\item In line 12 of Algorithm~\ref{alg2} we find the shortest path with the added  constraints Equations~(\ref{eq11}),(\ref{eq11a}) and $\alpha_{12} + \alpha_{25} \geq 1$ and $\alpha_{13} + \alpha_{35} \geq 1$. The path chosen is 1-3-4-5 of cost 4 where $\alpha_{35} =1$ and either $\alpha_{12}$ or $\alpha_{25}$ is also one.
\end{compactitem}
\noindent So the two best arcs to knockout 
so as to maximise the value of the shortest path after their elimination are arc (3,5) and either arc (1,2) or arc (2,5).

\FloatBarrier
 \clearpage
\newpage
 \pagestyle{empty}
\linespread{1}
\small \normalsize



\section*{Acknowledgments}
\noindent The author would like to acknowledge the comments made on an earlier version of this paper by anonymous reviewers.

\section*{Conflict of interest statement}
The author has no relevant financial or non-financial interests to disclose.


\end{document}